\documentclass[letterpaper, 10 pt, conference]{ieeeconf}
\IEEEoverridecommandlockouts
\usepackage[space]{cite}

\usepackage[tbtags]{amsmath}
\usepackage{amsmath}
\usepackage{amsfonts}
\usepackage{amssymb}
\usepackage{amsthm}
\usepackage{mathrsfs}
\usepackage{xspace}
\usepackage{xparse}
\usepackage{mathtools}

\usepackage{algorithm}
\usepackage{algpseudocode}
\usepackage{setspace}

\newtheorem{defn}{Definition}

\begin{document}

\title{
Data-Driven Stochastic Optimal Control Using Kernel Gradients
}

\author{
Adam~J.~Thorpe,~\IEEEmembership{Student~Member,~IEEE,}
Jake~A.~Gonzales,~\IEEEmembership{Student~Member,~IEEE,}\\
Meeko~M.~K.~Oishi,~\IEEEmembership{Senior Member,~IEEE}%
\thanks{%
    This material is based upon work supported by the National Science Foundation under NSF Grant Number CNS-1836900.  Any opinions, findings, and conclusions or recommendations expressed in this material are those of the authors and do not necessarily reflect the views of the National Science Foundation.
    The NASA University Leadership initiative (Grant \#80NSSC20M0163) provided funds to assist the authors with their research, but this article solely reflects the opinions and conclusions of its authors and not any NASA entity.
}
\thanks{A. Thorpe, J. Gonzales, and M. Oishi are with Electrical
  \& Computer Engineering, University of New Mexico, Albuquerque, NM. 
  Email: {\tt\{ajthor,jakegonzales,oishi\}@unm.edu}.
}
}

\maketitle

\begin{abstract}
We present an empirical, gradient-based method for solving data-driven stochastic optimal control problems using the theory of kernel embeddings of distributions. By embedding the integral operator of a stochastic kernel in a reproducing kernel Hilbert space, we can compute an empirical approximation of stochastic optimal control problems, which can then be solved efficiently using the properties of the RKHS. Existing approaches typically rely upon finite control spaces or optimize over policies with finite support to enable optimization. In contrast, our approach uses kernel-based gradients computed using observed data to approximate the cost surface of the optimal control problem, which can then be optimized using gradient descent. We apply our technique to the area of data-driven stochastic optimal control, and demonstrate our proposed approach on a linear regulation problem for comparison and on a nonlinear target tracking problem.
\end{abstract}


\section{Introduction}

The advent of autonomous systems, and the increasing complexity of real-world autonomy stemming from human interactions and learning-enabled components, obviates the need for algorithms which can accommodate real-world stochasticity. 
In such scenarios, model-based approaches may simply fail or hinge upon unrealistic assumptions such as linearity or Gaussianity, which can lead to questionable outcomes or unpredictable behaviors. 
One approach to dealing with such systems is data-driven control, which has proven to be useful for systems which may be resistant to traditional modeling techniques, or for which finding a simple mathematical model is simply impossible. 
In order to circumvent the problems faced by traditional model-based approaches, data-driven control uses empirical modeling techniques to synthesize implicit models which are amenable to analysis and control.
Nevertheless, these data-driven representations present new challenges for controller synthesis and optimization, which require the development of new tools and techniques to enable their use. 

We present a method for computing data-driven solutions to stochastic optimal control problems using an empirical, gradient-based approach. 
Our approach is based on Hilbert space embeddings of distributions, a nonparametric statistical learning technique that uses data collected from system observations to construct an implicit model of the dynamics as an element in a high-dimensional function space known as a reproducing kernel Hilbert space (RKHS).
Hilbert space embeddings of distributions have been applied to Markov models \cite{10.1145/1553374.1553497, 10.5555/3042573.3042778, 10.5555/3020652.3020720}, policy synthesis \cite{pmlr-v38-lever15, 10.5555/3060832.3060913}, state estimation and filtering \cite{10.5555/3104322.3104448, 6530747, 10.5555/2567709.2627677}, and also for solving stochastic optimal control problems \cite{thorpe2022stochastic, pmlr-v168-thorpe22a, nemmour2022maximum}. 
Additionally, these techniques admit finite sample bounds, which show convergence in probability as the sample size increases \cite{li2022optimal}.
Reproducing kernel Hilbert spaces and kernel embeddings of distributions, specifically, are broadly used in the area of nonparametric statistical inference and estimation. However, these techniques have not yet seen widespread adoption for controls.

The use of kernel methods for control has been explored in literature, and is closely related to the theory of Gaussian processes, Koopman operators, and support vector machines, in that they rely upon kernels or operators in high-dimensional function spaces.
The use of functional gradients in an RKHS for motion planning and policy synthesis have been used in \cite{Marinho-RSS-16, pmlr-v38-lever15}. However, these techniques either impose a particular problem structure, which limits their use more broadly, or rely upon a specific policy representation to compute the functional gradient.
Methods applying kernel embeddings of distributions to optimal control problems have been explored previously in \cite{nemmour2022maximum, 10.5555/3042573.3042778, pmlr-v168-thorpe22a} for MDPs and chance-constrained control. In addition, \cite{thorpe2022stochastic, pmlr-v168-thorpe22a} show that a kernel-based approximation of the stochastic optimal control problem can be solved as a linear program, but relies upon finite or discrete control spaces, which may be restrictive in practical control scenarios.

Our main contribution is a technique for computing solutions to stochastic optimal control problems using empirical, kernel-based stochastic gradient descent in an RKHS. 
Unlike existing functional gradient approaches such as \cite{pmlr-v38-lever15}, our approach does not rely upon explicit parameterizations of the policy in an RKHS.
Instead, we use the partial derivative reproducing property of kernels presented in \cite{zhou2008derivative} to compute an empirical gradient of the cost using observed data.
Our approach is based on the RKHS control framework presented in \cite{thorpe2022stochastic}, which optimizes over a finite set of user-specified admissible control actions. However, our proposed approach improves upon the techniques in \cite{thorpe2022stochastic} by eliminating the need for the control designer to strategically pre-select the policy support,
at the cost of increased computation time due to the iterative nature of gradient descent.

The rest of the paper is outlined as follows. 
In Section \ref{section: preliminaries}, we define the stochastic optimal control problem using kernel embeddings.
Then, in Section \ref{section: gradient method}, we describe the gradient-based optimization approach. 
In Section \ref{section: numerical results} we demonstrate our approach on a double integrator system for comparison to existing approaches and then on a nonholonomic vehicle system to demonstrate the capabilities of the approach. Concluding remarks are presented in Section \ref{section: conclusion}.


\section{Preliminaries \& Problem Formulation}
\label{section: preliminaries}

\begin{figure*}
    \centering
    \includegraphics[keepaspectratio,width=\textwidth]{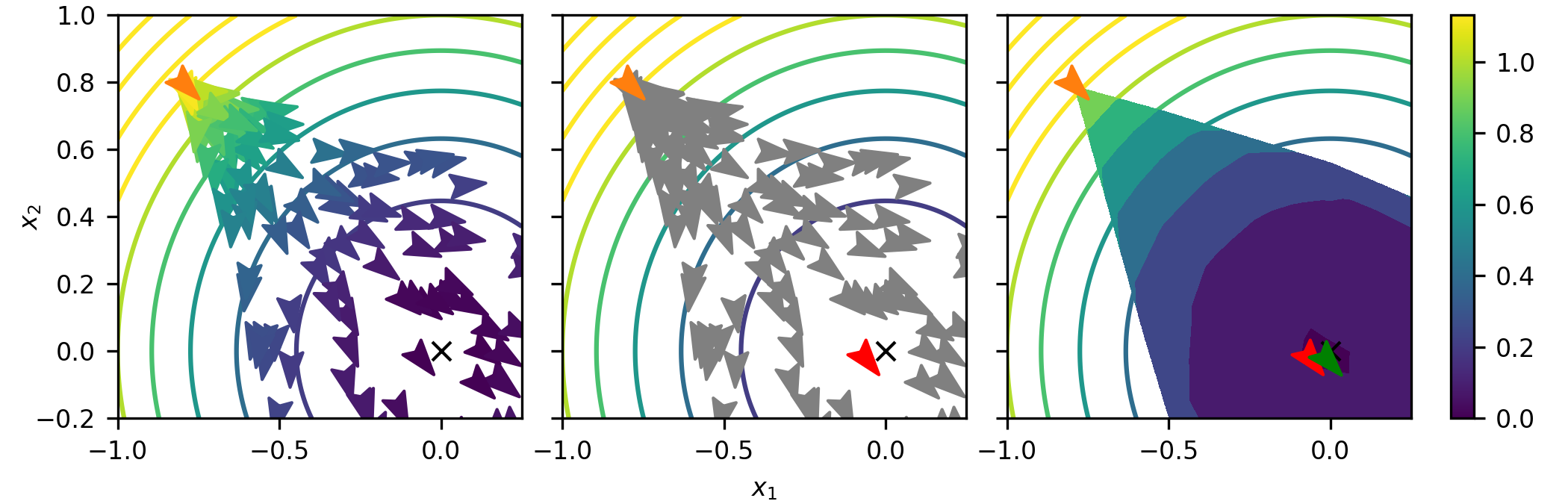}
    \caption{%
    Illustration of the gradient-based method on a stochastic optimal control problem of a nonholonomic vehicle system with bounded control authority seeking to minimize the Euclidean distance to the origin over a single time step. The initial condition is indicated by an orange arrow, the goal is denoted using an $\times$ at the origin, and the actual cost surface is depicted using contour lines. 
    (Left) Using the kernel-based estimate $\hat{m}$, we can empirically estimate the cost of taking control actions in an admissible set $\mathcal{A} \subset \mathcal{U}$. The resulting states after taking the actions in $\mathcal{A}$ from the initial condition are shown, color-coded by their estimated cost.
    (Center) The control algorithm in \cite{thorpe2022stochastic} chooses the control action in the admissible set $\mathcal{A}$ that minimizes the expected cost, but is sub-optimal. The resulting state after taking the chosen action is shown in red. 
    (Right) Our proposed approach using kernel-based gradient descent finds an approximately optimal solution by traversing the approximate cost surface (depicted using filled polygons), without resorting to a sampling-based approach. The resulting state after taking the approximately optimal control action is shown in green.
    }
    \label{fig: gradient big}
\end{figure*}


\subsection{System Model}

Let $(\mathcal{X}, \mathscr{B}_{\mathcal{X}})$ be a Borel space called the state space and $(\mathcal{U}, \mathscr{B}_{\mathcal{U}})$ be a compact Borel space called the control space.
Consider a discrete-time stochastic system,
\begin{equation}
    \label{eqn: system dynamics}
    x_{t+1} = f(x_{t}, u_{t}, w_{t}),
\end{equation}
where $x_{t} \in \mathcal{X} \subseteq \mathbb{R}^{n}$, $u_{t} \in \mathcal{U} \subset \mathbb{R}^{m}$, and $w_{t}$ are independent and identically distributed (i.i.d.) random variables representing a stochastic disturbance.
As shown in \cite{bertsekas1978stochastic}, the dynamics in \eqref{eqn: system dynamics} can equivalently be represented by a stochastic kernel $Q : \mathscr{B}_{\mathcal{X}} \times \mathcal{X} \times \mathcal{U} \to [0, 1]$ that assigns a probability measure $Q(\cdot \mid x, u)$ on $(\mathcal{X}, \mathscr{B}_{\mathcal{X}})$ to every $(x, u) \in \mathcal{X} \times \mathcal{U}$.

The system evolves from an initial condition $x_{0} \in \mathcal{X}$ (which may be drawn from an initial distribution $\mathbb{P}_{0}$ on $\mathcal{X}$) over a finite time horizon $t = 0, 1, \ldots, N$, $N \in \mathbb{N}$.


\subsection{Stochastic Optimal Control Problem}

Let $g : \mathcal{X} \to \mathbb{R}$ be an arbitrary \emph{convex} cost function, which we assume is measurable and bounded and lies in a Hilbert space of functions $\mathscr{H}$. 
At each time step, we seek the control $u \in \mathcal{U}$ that minimizes the following objective, 
\begin{align}
    \label{eqn: stochastic optimal control problem}
    \min_{u \in \mathcal{U}} \quad & \int_{\mathcal{X}} g(y) Q(\mathrm{d} y \mid x, u).
\end{align}

We assume that the stochastic kernel $Q$ is unknown, meaning we do not have direct information of the dynamics in \eqref{eqn: system dynamics} or the structure of the stochastic disturbance.
Instead, we assume that we have access to a sample $\mathcal{S} = \lbrace (x_{i}, u_{i}, y_{i}) \rbrace_{i=1}^{M}$ of observations taken i.i.d. from $Q$, 
where $x_{i}$ and $u_{i}$ are taken randomly from $\mathcal{X}$ and $\mathcal{U}$ and $y_{i} \sim Q(\cdot \mid x_{i}, u_{i})$.

Because the stochastic kernel $Q$ is unknown, we cannot solve \eqref{eqn: stochastic optimal control problem} directly since the integral in \eqref{eqn: stochastic optimal control problem} is intractable. 
Instead, as shown in \cite{thorpe2022stochastic}, we can use $\mathcal{S}$ to approximate the intractable integral with respect to $Q$ in \eqref{eqn: stochastic optimal control problem} as an empirical embedding in a high-dimensional space of functions known as a reproducing kernel Hilbert space. Then, we can solve an approximation of the original problem in \eqref{eqn: stochastic optimal control problem} in order to compute an approximately optimal control.
We outline the procedure below, but refer the reader to \cite{thorpe2022stochastic} for more details. 


\subsection{Approximate Problem Using Kernel Embeddings}

Define a positive definite kernel function $k : \mathcal{X} \times \mathcal{X} \to \mathbb{R}$ \cite[Definition~4.15]{steinwart2008support}. According to the Moore-Aronszajn theorem \cite{aronszajn1950theory}, given a positive definite kernel $k$, there exists a unique corresponding reproducing kernel Hilbert space (RKHS) $\mathscr{H}$ with $k$ as its reproducing kernel.

\begin{defn}
    A Hilbert space $\mathscr{H}$ of functions from $\mathcal{X}$ to $\mathbb{R}$ is called a reproducing kernel Hilbert space (RKHS) if there exists a positive definite function $k : \mathcal{X} \times \mathcal{X} \to \mathbb{R}$ called the reproducing kernel that satisfies the following properties:
    \begin{enumerate}
        \item 
        For every $x \in \mathcal{X}$, $k(x, \cdot) \in \mathscr{H}$, and 
        \item 
        For every $x \in \mathcal{X}$ and $f \in \mathscr{H}$, $f(x) = \langle f, k(x, \cdot) \rangle_{\mathscr{H}}$, which is known as the reproducing property.
    \end{enumerate}
\end{defn}

Similarly, we define the RKHS $\mathscr{U}$ of functions from $\mathcal{U}$ to $\mathbb{R}$ with $l : \mathcal{U} \times \mathcal{U} \to \mathbb{R}$ as its associated reproducing kernel.

According to \cite{10.1145/1553374.1553497, grunewalder2012conditional}, assuming the kernel $k$ is measurable and bounded, and given a probability measure $Q(\cdot \mid x, u)$, then by the Riesz representation theorem there exists a corresponding element $m(x, u) \in \mathscr{H}$ called the kernel distribution embedding, such that by the reproducing property, $\langle g, m(x, u) \rangle_{\mathscr{H}} = \int_{\mathcal{X}} g(y) Q(\mathrm{d} y \mid x, u)$.
This means that by representing the integral operator with respect to $Q$ as an element in the RKHS, we can compute the expectation of any function $f \in \mathscr{H}$ as an RKHS inner product. 

Using a sample $\mathcal{S}$, we can compute an empirical estimate $\hat{m}(x, u)$ of $m(x, u)$ as the solution to a regularized least-squares problem \cite{grunewalder2012conditional}.
The solution is given by,
\begin{equation}
    \label{eqn: empirical embedding}
    \hat{m}(x, u) = \Phi^{\top} (G + \lambda M I)^{-1} \Psi k(x, \cdot) l(u, \cdot),
\end{equation}
where $\Phi, \Psi$ are feature vectors with elements $\Phi_{i} = k(y_{i}, \cdot)$ and $\Psi_{i} = k(x_{i}, \cdot) l(u_{i}, \cdot)$, respectively, and $G \in \mathbb{R}^{M \times M}$ is a real matrix where the $ij^{\rm th}$ element is $k(x_{i}, x_{j}) l(u_{i}, u_{j})$. 
For simplicity of notation, let $W = (G + \lambda M I)^{-1}$.

Using the estimate $\hat{m}(x, u)$, we can approximate the intractable integrals with respect to $Q$ in \eqref{eqn: stochastic optimal control problem} via an RKHS inner product, 
\begin{align}
    \int_{\mathcal{X}} g(y) Q(\mathrm{d} y \mid x, u) &\approx \langle g, \hat{m}(x, u) \rangle_{\mathscr{H}} \\
    \label{eqn: kernel approximation}
    &= \boldsymbol{g}^{\top} W \Psi k(x, \cdot) l(u, \cdot),
\end{align}
where $\boldsymbol{g}$ is a vector with elements $\boldsymbol{g}_{i} = g(y_{i})$.

This representation is key to our approach, since it means we can approximate the previously intractable problem in \eqref{eqn: stochastic optimal control problem} using data comprised of system observations. 


\subsection{Problem Statement}

Following \cite{thorpe2022stochastic}, we can approximate the stochastic optimal control problem \eqref{eqn: stochastic optimal control problem} using the estimate $\hat{m}(x, u)$ and \eqref{eqn: kernel approximation} as,
\begin{align}
    \label{eqn: approximate stochastic optimal control problem}
    \min_{u \in \mathcal{U}} \quad & \boldsymbol{g}^{\top} W \Psi k(x, \cdot) l(u, \cdot).
\end{align}
Theoretically, we could optimize for $u$ directly. 
However, this is a non-convex problem in general, which makes solving \eqref{eqn: approximate stochastic optimal control problem} difficult. 
For example, it may be exceptionally difficult to solve \eqref{eqn: approximate stochastic optimal control problem} for common kernel choices such as the Gaussian kernel function, since optimizing a linear combination of Gaussians is a non-convex problem. 
Thus, finding a control action $u \in \mathcal{U}$ that minimizes the approximate problem presents a significant challenge.

One possible approach is given in \cite{thorpe2022stochastic, pmlr-v168-thorpe22a}, where a stochastic control policy with finite support is represented as an embedding in the RKHS $\mathscr{U}$. Then, the policy can be obtained as the solution to a linear program, giving a set of probability values over a user-specified set of admissible control actions $\mathcal{A} \subset \mathcal{U}$.
However, a significant drawback of this approach is the need to strategically select $\mathcal{A}$ such that it contains controls which are close to the true solution, and typically only finds a sub-optimal solution to the approximate problem. 

We propose to compute the control input via a kernel-based stochastic gradient descent method.
Using the properties of reproducing kernel Hilbert spaces, we can compute the gradient by taking the partial derivative of the kernel, rather than explicitly computing the gradient with respect to $u$. This allows us to optimize the control input by directly optimizing within the RKHS and avoids the problem of non-convexity in optimizing for $u$ in the approximate problem. 
An illustration of this idea is depicted in Figure \ref{fig: gradient big}.


\section{Computing Controls Using \\ Gradient Descent in an RKHS}
\label{section: gradient method}

We seek to compute the partial derivative of the objective in \eqref{eqn: approximate stochastic optimal control problem} with respect to $u$. 
We first define the notation used to describe the partial derivative of a bivariate function.

\begin{defn}[Partial Functional Derivative Notation]
    Given a bivariate function $l : \mathcal{U} \times \mathcal{U} \to \mathbb{R}$, $u, u' \in \mathcal{U}$, we denote the partial derivative as,
    \begin{equation}
        \partial^{p, q} l(u, u') \coloneqq \frac{\partial^{p_{1} + \cdots + p_{m} + q_{1} + \cdots + q_{m}} l(u, u')}{\partial(u_{1})^{p_{1}} \cdots \partial(u_{m})^{p_{m}} \partial(u_{1}^{\prime})^{q_{1}} \cdots \partial(u_{m}^{\prime})^{q_{m}}},
    \end{equation}
    where $p, q$ are multi-indices.
\end{defn}

As shown in \cite{zhou2008derivative}, we can compute the partial derivative of any function $h \in \mathscr{U}$ via the reproducing property as,
\begin{equation}
    \label{eqn: partial derivative reproducing property}
    \partial^{p} h(u) = \langle h, \partial^{p, 0} l(u, \cdot) \rangle_{\mathscr{U}}.
\end{equation}
In short, this means that we do not need to directly compute the partial derivative of the cost function $g$ with respect to $u$ (which may be unknown if we are only given points $g(y_{i})$), and we may compute the empirical gradient by taking the partial derivative of the kernel $l$. 

Let $\hat{J}(u) = \boldsymbol{g}^{\top} W \Psi k(x, \cdot) l(u, \cdot)$ be the objective of the approximate optimal control problem in \eqref{eqn: approximate stochastic optimal control problem}.
Note that $\hat{J}(u)$ can be written using the reproducing property as 
\begin{equation}
    \hat{J}(u) 
    = \langle \boldsymbol{g}^{\top} W \Psi k(x, \cdot), l(u, \cdot) \rangle_{\mathscr{U}}.
\end{equation}
Then, using \eqref{eqn: partial derivative reproducing property}, the partial derivative of $\hat{J}(u)$ with respect to the control $u$ can be computed as,
\begin{equation}
    \label{eqn: empirical cost gradient}
    \partial^{1} \hat{J}(u) = \langle \boldsymbol{g}^{\top} W \Psi k(x, \cdot), \partial^{1, 0} l(u, \cdot) \rangle_{\mathscr{U}}.
\end{equation}

This approach has a significant advantage, most notably that most popular kernels are easy to differentiate, meaning we can quickly compute the empirical gradient for an arbitrary cost function $g \in \mathscr{H}$. In addition, the empirical cost gradient can be computed as a simple matrix multiplication. 

As a practical example, consider the Gaussian kernel, $l(u, u') = \exp(-\lVert u - u' \rVert_{2}^{2}/ 2 \sigma^{2})$, $\sigma > 0$ (assuming $u$ is a scalar variable for simplicity). The partial derivative of the Gaussian kernel is given by, 
\begin{equation}
    \partial^{1, 0} l(u, u') = - \frac{\lvert u - u' \rvert}{\sigma^{2}} \exp \biggl( -\frac{\lVert u - u' \rVert^{2}}{2 \sigma^{2}} \biggr).
\end{equation}
Then the partial derivative of the objective in \eqref{eqn: empirical cost gradient} can be computed as $\partial^{1} \hat{J}(u) = \boldsymbol{g}^{\top} W (\Psi k(x, \cdot) l(u, \cdot) \odot \Delta)$,
where $\odot$ denotes the Hadamard (or element-wise) product and $\Delta \in \mathbb{R}^{M}$ is a vector with elements $\Delta_{i} = - \lvert u_{i} - u \rvert / \sigma^{2}$. 

We use the empirical gradient of the cost function $g$ computed using \eqref{eqn: empirical cost gradient} in order to compute the gradient direction for stochastic gradient descent. Then, by traversing the approximate cost surface using the empirical gradient, we obtain an approximately optimal solution to the problem in \eqref{eqn: approximate stochastic optimal control problem}.
We outline the procedure in Algorithm \ref{alg: kernel gradient descent}.

\begin{algorithm}
\caption{Kernel-Based Gradient Descent}
\label{alg: kernel gradient descent}
\begin{algorithmic}[1]
    \State \textbf{given} 
    embedding estimate $\hat{m}$,
    initial guess $u_{0}$
    \Repeat
        \State $\Delta u_{n} \gets \langle  \boldsymbol{g}^{\top} W \Psi k(x, \cdot), \partial^{1, 0} l(u_{n}, \cdot) \rangle_{\mathscr{U}}$
        \State \textbf{choose} step size $\eta$
        \State $u_{n+1} \gets u_{n} - \eta \Delta u_{n}$
    \Until{stopping criterion satisfied}
    \State \textbf{return} $u_{n}$
\end{algorithmic}
\end{algorithm}

Since the estimate $\hat{m}$ converges in probability to the true embedding $m$ at a minimax optimal rate of $\mathcal{O}(M^{-1/2})$ \cite{li2022optimal}, the approximate cost surface also converges in probability to the true cost surface. Hence, as the sample size increases, we obtain a closer approximation of the true cost surface. 
However, it is important to note that the empirical cost surface is generally not convex, even if the original function is convex, meaning we are only guaranteed to find a locally optimal solution to the approximate problem. This is obvious, since the noise of the data also adds noise to the empirical cost surface. Nevertheless, we can use more advanced gradient descent methods (e.g. using momentum or a ``temperature'' in place of the learning rate) to mitigate the issues of optimizing over an empirical cost surface. This also motivates the need to choose an initial guess as close as possible to the optimal solution, which is detailed in the next section. 


\subsection{Initialization}

Initializing the gradient descent algorithm close to the true solution ensures that we obtain an approximately optimal solution in fewer gradient steps. One possibility is to compute a sub-optimal initial guess for Algorithm \ref{alg: kernel gradient descent} using \cite{thorpe2022stochastic}.

As shown in \cite{thorpe2022stochastic}, we can compute a solution to the (unconstrained) approximate stochastic optimal control problem in \eqref{eqn: approximate stochastic optimal control problem} by representing a stochastic policy $\pi : \mathscr{B}_{\mathcal{U}} \times \mathcal{X} \to [0, 1]$ as a kernel embedding $p$ in the RKHS $\mathscr{U}$, 
\begin{equation}
    \label{eqn: kernel stochastic policy}
    p(x) = \sum_{j=1}^{P} \gamma_{j}(x) l(\tilde{u}_{j}, \cdot) = \Upsilon^{\top} \gamma(x),
\end{equation}
where $\gamma(x) \in \mathbb{R}^{P}$ are real-valued coefficients that depend on the state $x \in \mathcal{X}$, $\Upsilon$ is a feature vector with elements $\Upsilon_{j} = l(\tilde{u}_{j}, \cdot)$, and the points $\mathcal{A} = \lbrace \tilde{u}_{j} \rbrace_{j=1}^{P}$ are a set of user-specified admissible control actions that we want to optimize over. The problem then becomes finding the coefficients $\gamma(x)$ that optimize the approximate control problem. According to \cite{thorpe2022stochastic}, we can view the coefficients $\gamma(x)$ as a set of probabilities that weight the user-specified control actions in $\mathcal{A}$, which we can find as the solution to a linear program,
\begin{subequations}
\label{eqn: linear program}
\begin{align}
    \min_{\gamma(x) \in \mathbb{R}^{P}} \quad & \boldsymbol{g}^{\top} W \Psi k(x, \cdot) \Upsilon^{\top} \gamma(x) \\
    \text{s.t.} \quad & \boldsymbol{1}^{\top} \gamma(x) = 1 \\
    & 0 \preceq \gamma(x)
\end{align}
\end{subequations}
The linear program can efficiently be solved via the Lagrangian dual. Letting $C(x) = \boldsymbol{g}^{\top} W \Psi k(x, \cdot) \Upsilon^{\top}$, the solution according to \cite{boyd2004convex} is given by a vector of all zeros except at the index $j = \arg\min_{i} \lbrace C_{i}(x)^{\top} \rbrace$, where it is $1$. In other words, we choose the control action in $\mathcal{A}$ that corresponds to the index $j$ which is the solution to the Lagrangian dual problem.
See \cite{thorpe2022stochastic} for more details.

By choosing control actions in $\mathcal{A}$ that are good candidate solutions to the optimal control problem in \eqref{eqn: stochastic optimal control problem}, we obtain a good initial guess for the gradient-based learning algorithm. However, unlike the approach in \cite{thorpe2022stochastic}, we do not require the approximately optimal solution to lie within $\mathcal{A}$, and we further improve the solution of the LP using stochastic gradient descent. 


\section{Numerical Results}
\label{section: numerical results}

We demonstrate our approach on a regulation problem using a discrete-time stochastic chain of integrators for verification, and on a target tracking problem using nonholonomic vehicle dynamics to demonstrate the utility of the approach. 
For all problems, we use a Gaussian kernel for $k$ and $l$, which has the form $k(x, x') = \exp( -\lVert x - x' \rVert_{2}^{2}/2 \sigma^{2})$,
where $\sigma > 0$.
Following \cite{10.1145/1553374.1553497}, we choose the regularization parameter to be $\lambda = 1/M^{2}$, where $M \in \mathbb{N}$ is the sample size used to construct the estimate $\hat{m}$. In practice, the parameters $\sigma$ and $\lambda$ are typically chosen via cross-validation, where $\sigma$ is chosen according to the relative spacing of the data points (usually the median distance) and $\lambda$ is chosen such that $\lambda \to 0$ as $M \to \infty$.
A more detailed discussion of parameter selection is outside the scope of the current work (see \cite{li2022optimal} for recent results on regularization rates). 
Numerical experiments were performed in Python on an AWS cloud computing instance.
Code for all analysis and experiments is available as part of the stochastic optimal control using kernel methods (SOCKS) toolbox \cite{10.1145/3501710.3519525}. 


\subsection{Regulation of a Double Integrator System}

We consider the problem of regulation for a 2D stochastic chain of integrators system, with dynamics given by 
\begin{equation}
    x_{t+1} = 
    \begin{bmatrix}
        1 & T_{s} \\
        0 & 1
    \end{bmatrix}
    x_{t} + 
    \begin{bmatrix}
        T_{s}^{2}/2 \\
        T_{s}
    \end{bmatrix}
    u_{t} + w_{t},
\end{equation}
where $x_{t} \in \mathbb{R}^{2}$ is the state, $u_{t} \in \mathbb{R}$ is the control input, which we constrain to be within $[-1, 1]$, $w_{t}$ is a random variable with distribution $\mathcal{N}(0, 0.01 I)$, and $T_{s}$ is the sampling time.  
We seek to compute a control input $u$ as the solution to the following stochastic optimal control problem,
\begin{equation}
    \min_{u} \quad \int_{\mathcal{X}} g(y) Q(\mathrm{d} y \mid x, u),
\end{equation}
where $Q$ is a representation of the dynamics as a stochastic kernel.
We use the cost function $g(x) = \lVert x \rVert_{2}$, which serves to drive the system to the origin. 

We consider a sample $\mathcal{S} = \lbrace (x_{i}, u_{i}, y_{i}) \rbrace_{i=1}^{M}$ of size $M = 1600$ taken i.i.d. from $Q$. The states $x_{i}$ were taken uniformly in the region $[-1, 1] \times [-1, 1]$, the control inputs $u_{i}$ were taken uniformly from $[-1, 1]$, and the resulting states were generated according to $y_{i} \sim Q(\cdot \mid x_{i}, u_{i})$.
We then presumed no knowledge of the system dynamics or the stochastic disturbance for the purpose of computing the approximately optimal control action using our proposed method. 
Using the sample $\mathcal{S}$, we then computed an estimate $\hat{m}$ of the kernel embedding $m$ as in \eqref{eqn: empirical embedding} using Gaussian kernels with bandwidth parameter $\sigma = 3$, chosen via cross-validation. We then selected $R = 121$ \emph{evaluation points} $\lbrace x_{j} \rbrace_{j=1}^{R}$ spaced uniformly in the region $[-1,1] \times [-1,1]$, from which to compute the optimal control actions. 

\begin{figure}
    \centering
    \includegraphics[keepaspectratio,width=\columnwidth]{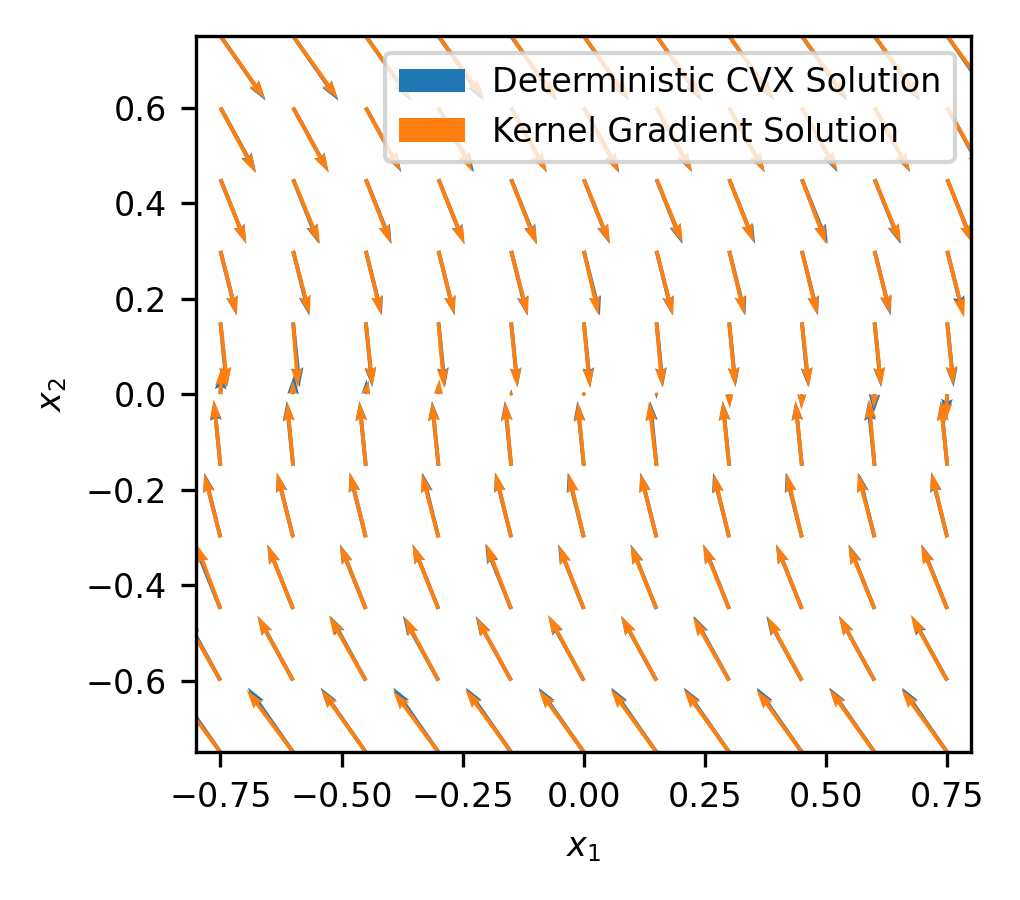}
    \caption{Vector field of the closed-loop dynamics of a deterministic double integrator system under an optimal control strategy computed using CVX (blue)
    and the vector field of the approximately optimal closed-loop system under the gradient descent-based control algorithm (orange). We can see that the kernel-based gradient descent solution closely matches the solution from CVX.}
    \label{fig: vector field plot}
\end{figure}

\begin{figure}
    \centering
    \includegraphics[keepaspectratio,width=\columnwidth,height=2in]{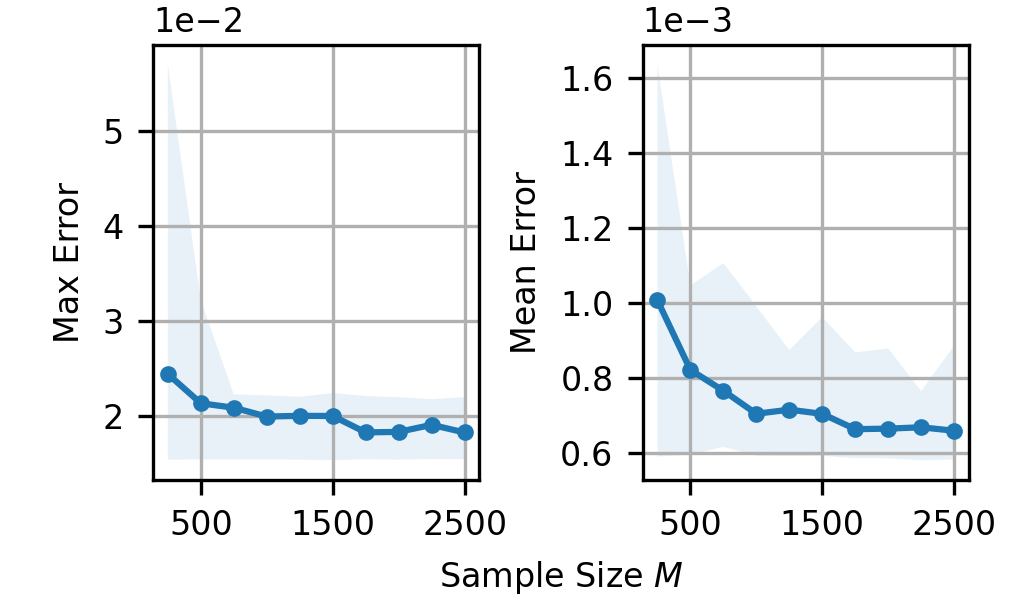}
    \caption{(Left) Maximum error of the control inputs computed via our proposed method versus the optimal control inputs computed via CVX for varying sample sizes $M \in [250, 2500]$. (Right) The mean error of the control inputs.}
    \label{fig: error plot}
\end{figure}

To provide a basis for comparison, we computed the optimal control actions using CVX from the evaluation points $\lbrace x_j \rbrace_{j=1}^{R}$ using the \emph{deterministic} dynamics. We then propagate the dynamics forward in time using the optimal inputs to obtain the state at the next time instant. The vector field of the closed-loop dynamics under the optimal control inputs is shown in Figure \ref{fig: vector field plot} (blue). 
We then computed the approximately optimal control actions using Algorithm \ref{alg: kernel gradient descent} with the sample $\mathcal{S}$ taken from the \emph{stochastic} dynamics to minimize the cost at each point $\lbrace x_j \rbrace_{j=1}^{R}$ over a single time step. For Algorithm \ref{alg: kernel gradient descent}, we used a step size of $\eta=0.01$ and limited the number of iterations to $100$. 
The vector field of the closed-loop dynamics using the approximately optimal solution computed using our proposed method is shown in Figure \ref{fig: vector field plot} (orange). 

We can see that the gradient based algorithm computes approximately optimal control inputs which closely match the solution computed via CVX. 
This demonstrates the effectiveness of the gradient-based algorithm to compute approximately optimal control actions with no prior knowledge of the dynamics or the stochastic disturbance. 

Note that the quality of the empirical approximation of the cost surface depends on the sample size $M$. As the sample size increases, the approximation improves, and we obtain a closer approximation of the optimal solution using our method. To demonstrate this, we computed the mean error and the maximum error between the \emph{data-driven} gradient-based solution and the solution via CVX using the \emph{deterministic} dynamics for varying sample sizes $M \in [250, 2500]$ averaged over $20$ iterations. The results are shown in Figure \ref{fig: error plot}. We can see that the error of the approximately optimal solution decreases as $M$ increases. However, we can also see that the quality of the solution does not improve appreciably as the sample size increases, which is due to the asymptotic convergence of the estimate $\hat{m}$ to the true embedding $m$. This presents a tradeoff between computation time and numerical accuracy, since the computational complexity scales polynomially with the sample size.


\subsection{Target Tracking Using a Nonholonomic Vehicle}

\begin{figure}
    \centering
    \includegraphics[keepaspectratio,width=\columnwidth]{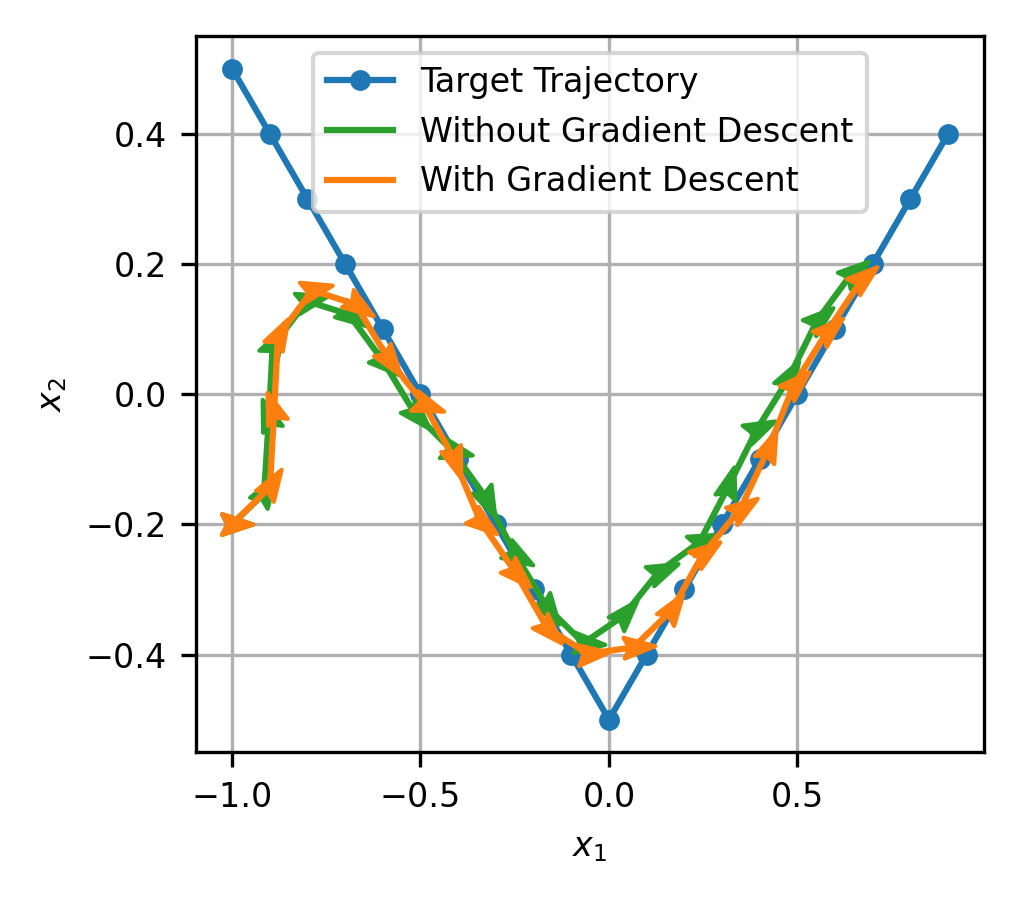}
    \caption{Trajectory generated via our proposed method (orange) which tracks the target trajectory (blue). The trajectory generated using the technique in \cite{thorpe2022stochastic} is shown for comparison (green). Note that the trajectory computed using our proposed approach more closely follows the target trajectory.}
    \label{fig: nonholonomic example}
\end{figure}

We consider the problem of target tracking for a nonholonomic vehicle system as in \cite{thorpe2022stochastic}. The dynamics are given by
\begin{align}
    \dot{x}_{1} = u_{1} \sin(x_{3}), &&
    \dot{x}_{2} = u_{1} \cos(x_{3}), &&
    \dot{x}_{3} = u_{2}
\end{align}
where $x = [x_{1}, x_{2}, x_{3}]^{\top} \in \mathbb{R}^{3}$ are the states, $u = [u_{1}, u_{2}]^{\top} \in \mathbb{R}^{2}$ are the control inputs. 
The control inputs are constrained such that $u_{t} \in [0.5, 1.2] \times [-10.1, 10.1]$.
We discretize the dynamics in time using a zero-order input hold and apply an affine disturbance $w \sim \mathcal{N}(0, 0.01 I)$. 
We define a target trajectory as a sequence of position coordinates indexed by time, shown in Figure \ref{fig: nonholonomic example} (blue).
We choose an initial condition of $x_{0} = [-1, -0.2, \pi/2]^{\top}$, and evolve the system forward in time over the time horizon $N = 20$.
At each time $t$, starting at $t = 0$, we seek to compute a control input $u_{t}$ as the solution to the following stochastic optimal control problem,
\begin{align}
    \min_{u_{t}} \quad \int_{\mathcal{X}} g_{t}(y) Q(\mathrm{d} y \mid x_{t}, u_{t}),
\end{align}
where the cost function $g_{t}(x)$ seeks to minimize the squared Euclidean distance of the system's position to the target trajectory's position at each time step. 

We consider a sample $\mathcal{S} = \lbrace (x_{i}, u_{i}, y_{i}) \rbrace_{i=1}^{M}$, with sample size $M = 3000$. 
The states $x_{i}$ were taken uniformly within the region $[-1.2, 1.2] \times [-1.2, 1.2] \times [-\pi, \pi]$, the control actions $u_{i}$ were taken uniformly within the region $[0.5, 1.2] \times [-10.1, 10.1]$, and the resulting states were drawn according to $y_{i} \sim Q(\cdot \mid x_{i}, u_{i})$.
Using the sample $\mathcal{S}$, we then computed an estimate $\hat{m}$ of the kernel embedding $m$ as in \eqref{eqn: empirical embedding} using Gaussian kernels with bandwidth parameter $\sigma = 3$, chosen via cross-validation.

In order to compute a baseline for comparison, we then computed the control actions using \cite{thorpe2022stochastic}, which computes a stochastic policy embedding $p_{t}$ at every time step over a finite set $\mathcal{A} = \lbrace \tilde{u}_{j} \rbrace_{j=1}^{210}$, as in \eqref{eqn: kernel stochastic policy}. We chose the controls $\tilde{u}_{j}$ in the admissible set $\mathcal{A}$ to be uniformly spaced in the range $\tilde{u}_{j} \in [0.5, 1.2] \times [-10.1, 10.1]$. Starting at the initial condition $x_{0}$, we then computed the control actions by solving \eqref{eqn: linear program} at each time step forward in time. The resulting trajectory is plotted in Figure \ref{fig: nonholonomic example} (green), had a total cost of $\sum_{t=0}^{N} g_{t}(x_{t}) = 1.447$, and the computation time was $1.159$ seconds.

We then evolve the system forward in time using the approximately optimal control action selected via the kernel-based gradient descent algorithm (Algorithm \ref{alg: kernel gradient descent}), and initializing using the solution to \eqref{eqn: linear program} using $\mathcal{A}$ as above. For Algorithm \ref{alg: kernel gradient descent}, we chose a step size of $\eta = 0.1$ and limited the number of gradient iterations to $100$. The resulting trajectory is plotted in Figure \ref{fig: nonholonomic example} (orange), and has a total cost of $\sum_{t=0}^{N} g_{t}(x_{t}) = 1.363$. The total computation time was $8.865$ seconds. As expected, we see that the trajectory computed using our method more closely follows the target trajectory (has a lower overall cost). This shows that the gradient-based algorithm is able to compute the approximately optimal control actions for a nonlinear system at each time step, using only data collected from system observations. 


\section{Conclusions \& Future Work}
\label{section: conclusion}

In this paper, we presented a method for computing the approximately optimal control action for stochastic optimal control problems using a data-driven approach. Our proposed method leverages kernel gradient-based methods and achieves more optimal control solutions than existing sample-based approaches. We plan to explore methods to compute control solutions more efficiently using the geometric properties of the RKHS, e.g. via projections, and to adapt the algorithm to dynamic programs and constrained stochastic optimal control problems. 


\bibliographystyle{IEEEtran}
\bibliography{IEEEabrv, shortIEEE, bibliography}

\end{document}